# Use of spheroidal models in gravitational tomography

Valery Sizikov

Saint-Petersburg University ITMO, Saint-Petersburg, Russia
Email: sizikov2000@mail.ru

**Abstract**

The direct gravimetry problem is solved using the subdivision of each body of a deposit into a set of vertical adjoining bars, and in the inverse problem each body of a deposit is modeled by a uniform ellipsoid of revolution (spheroid). Well-known formulas for $z$-component of gravitational intensity of a spheroid are transformed to a convenient form. Parameters of a spheroid are determined by minimizing the Tikhonov smoothing functional using constraints on the parameters. This makes the ill-posed inverse problem by unique and stable. The Bulakh algorithm for initial estimating the depth and mass of a deposit is modified. The technique is illustrated by numerical model examples of deposits in the form of two and five bodies. The inverse gravimetry problem is interpreted as a gravitational tomography problem or the intravision of the Earth's crust and mantle.

**Keywords:**

direct and inverse gravimetry problems, modeling of deposit bodies by spheroids, gravitational tomography, Earth's intravision

## 1. Introduction

Modeling of deposits is one of the main ways for solving the direct and especially inverse gravimetry problems. The direct problem is a calculation of the gravitational field produced by some model deposit on the Earth's surface (when performing practical measurements of the field, the direct problem is not solved). The inverse (more complicated) problem is a determination by mathematical and computer way of the deposit parameters on the basis of calculated or measured field anomaly (e.g., the Bouguer's anomaly) on the Earth's surface. In solving the direct problem, it is desirable to approximate a deposit by a few bodies of an enough arbitrary shape. In solving the inverse problem, one should use bodies of a more or less regular form close to the deposit.

### (a) Models of deposits

Various authors use the following *models of deposits* [1–3]: in the form of quadrangular truncated pyramids, prisms, cylinders, beams, polyhedrons, parallelepipeds, intersecting bars, etc. However, such figures have nonsmooth surfaces and generate cumbersome (although not complicated) formulas (see, e.g., [1]). We note also a plane-layered model [4]. In the works [5–8], the modeling of deposits uses homogeneous (and inhomogeneous) *spheroids*, or ellipsoids of revolution which are effectively applied, for example, in astrophysics for constructing galactic models [9, 10]. In this paper, we continue to use spheroids for deposit modeling.

To calculate a deposit, a different information is usually used, namely, intensities of the gravitational and magnetic fields [1, 2, 5–8, 11–13], gravity gradient tensor components [3], seismic data [11, 14], remote sensing from satellites [11, 13], et al. We will use only vertical component of the gravitational field intensity $V_z$ on the Earth's surface and this will be, in principle, sufficient for calculating enough complicated deposit models.

In [14–18] et al., variants for determining a boundary shape $z(x)$ separating two parts of the Earth's crust (properly crust and deposit) are stated. Experimental function is an anomaly of the gravitational force $\Delta g(x)$ (e.g., the Bouguer's anomaly [13]). The lower boundary of the deposit $H = \text{const}$ [15, 16] or the coordinates of its center [17, 18], as well as the density anomaly $\Delta\rho$ are additionally prescribed. One solves an one-dimensional nonlinear integral equation for $z(x)$, where $z(x)$ is an upper bound of the deposit [15, 16] or of the entire boundary [17, 18]. Furthermore, each cross-section of the deposit is modeled by an ellipse in [17] or the boundary is arbitrary in [18]. In [15–18], the problem is solved as a set of one-dimensional problems (in a number of vertical cross-sections). It is important that a number of parameters are a priori prescribed, namely, the lower boundary of the body *H*, the coordinates of its center and the density difference $\Delta\rho$. As a result, the technique for solving the inverse problem becomes limited, although one-valued. We desire to solve the inverse problem, not prescribing, but determining the deposit parameters at the expense of decremental constraints.

In solving the *direct problem* (calculation of the gravity field anomaly produced by the deposit on the Earth's surface), we will consider the deposit in the form of a few bodies and subdivide each body into a set of adjoining *vertical bars* [7]. And in solving the *inverse problem* (determination of the deposit parameters from the field anomaly measured at the Earth's surface), we will simulate each body of the deposit by a *biaxial ellipsoid*, or *ellipsoid of revolution*, or *spheroid* (the convenient astrophysical term).

Spheroids (ellipsoids) are widely used in celestial mechanics [19, 20], astrophysics (galactic models) [9, 10], and geophysics [21]. However, in many works, formulas for the field of an ellipsoid are not reduced to a convenient form. Moreover, after the works of Yun'kov [21], spheroids have been used in geophysics not often [6, 8, 22, 23].

In [22], it is assumed that bodies of the ore type which are sources of gravitational field have the form closed to spheroids. Spheroids are homogeneous, convex and star-shaped domains possessing the mean plane and uniqueness theorems hold for such domains. In the inverse problem, forms of the bodies are determined via minimizing the discrepancy functional using the Lagrange undetermined multipliers (the regularization parameters). However, coordinates of the centers and densities of the bodies are prescribed, which limits the problem.

In [23], it is assumed that local inclusions have the form of homogeneous bodies of revolution, in particular, sphere, oblate or prolate spheroid, and the solution of the inverse problem is sought in the form of a series in polynomials and associated Legendre functions, as well as in the form of splines via minimizing the Tikhonov functional by a variation method.

In the direct problem, we will assume that geologic bodies have an enough arbitrary shape (although resembling a spheroid). And in the inverse problem, we will model them by spheroids. Furthermore, we will include into a number of unknowns coordinates of the centers of spheroids, their semiaxis and densities (with constraints), as well as their quantity.

### (b) Comparison with tomography

The inverse gravimetry problem is often solved as a set of two-dimensional problems, namely, the density anomaly distribution in a number of vertical cross-sections are determined and then a three-dimensional (volume) picture is composed. This procedure reminds the techniques that are typical for different types of tomography [24–28], first of all, X-ray computerized tomography (XCT) and NMR tomography (MRT) [24–27]. In this paper (as in [6–8]), the spheroids are used and three-dimensional problem is solved and this is equivalent to the three-dimensional tomography [29]. Therefore, it is proposed to refer the inverse gravimetry problem to the tomography problems (as this has already done in [6–8]) and to call the inverse gravimetry problem as the *gravitational tomography problem*. Moreover, the peculiarity of the given problem is that it gives the possibility to 'look' into the Earth's depth via mathematical processing the surface results (without drilling holes), and this reminds the operation of intravision which is also typical for the tomography. Therefore, it is proposed to call the inverse gravimetry problem also the *problem of intravising the Earth* [6–8]. The proposed interpretations of the inverse gravimetry problem will allow to use extensive developments in the field of computerized tomography, in particular, in XCT and MRT.

## 2. Calculating the direct problem with use of vertical bars

Consider a geologic deposit in the form of several homogeneous bodies having an enough arbitrary form. In figure 1, as an example, we adduce a deposit in the form of two bodies. We will conditionally associate body 1 with an ore body and body 2 with an intrusion.

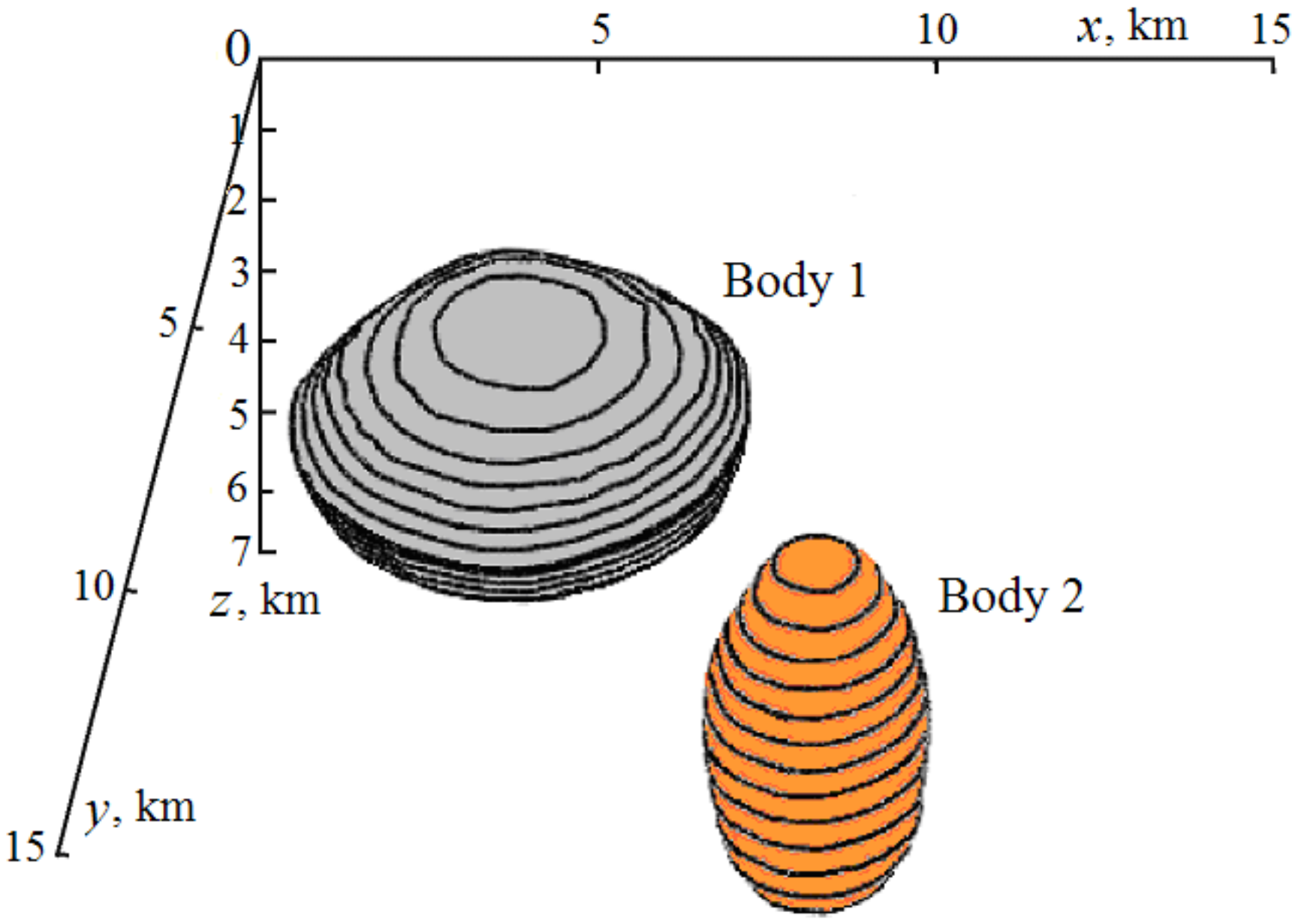

**Figure 1**. Deposit in the form of two bodies. (Online version in colour.)

In figure 2, the contours of $z$-sections of each body are shown.

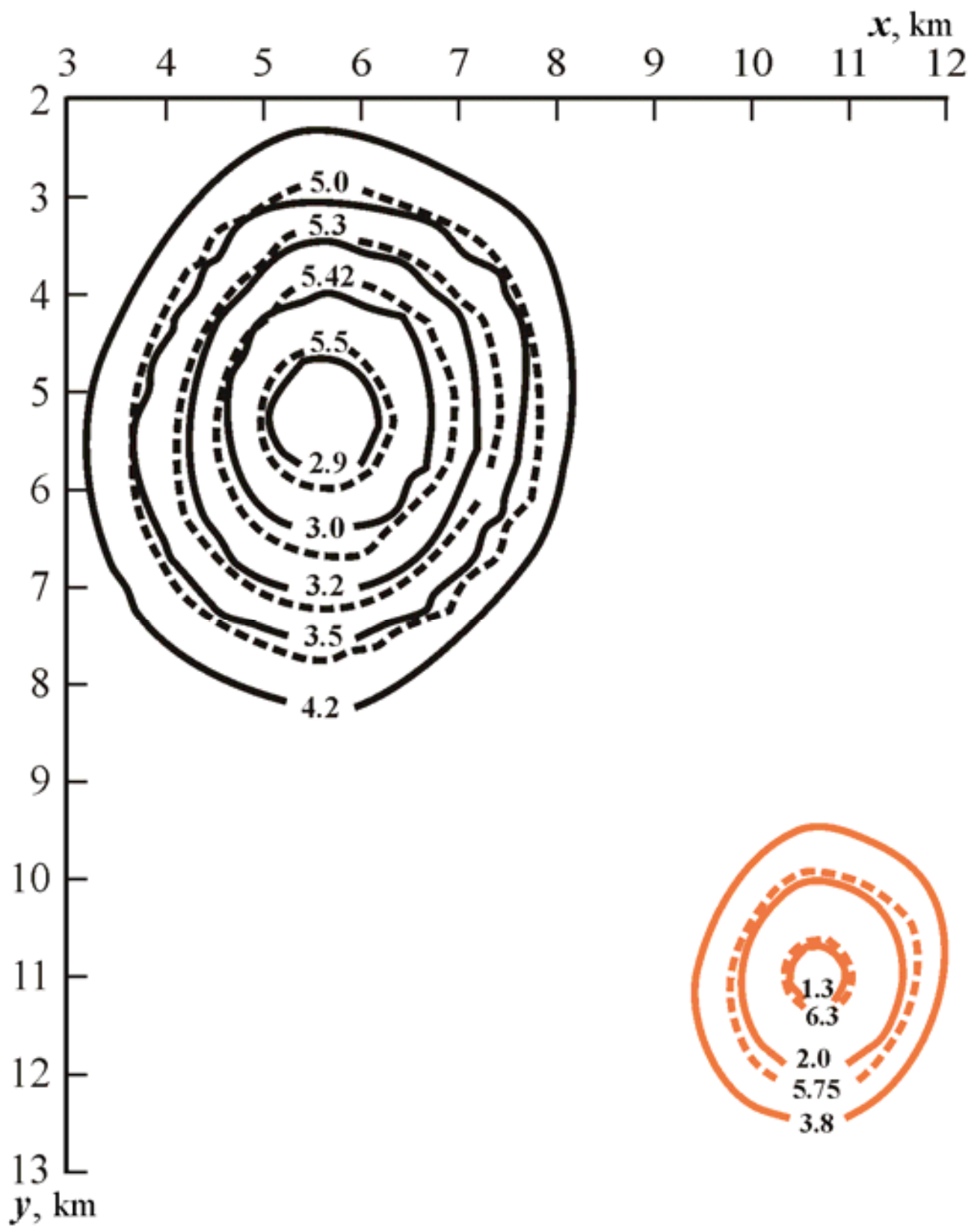

**Figure 2**. Contours of $z$-sections of the deposit bodies. (Online version in colour.)

Numbers on the contours (lines) are *z*-coordinates (in km) of the sections. Solid lines are contours located above the conditional median section (mean plane via the Sretenskii class of geologic models [30, 31]), and dashed lines are contours under one.

**Definition 2.1** [7, 8]. A body is called *vertically star-shaped* if any vertical ray (straight line) intersects its boundary only twice.

Let the density of a body $\rho = \text{const}$ and the body be vertically star-shaped. We represent it by a collection of vertical elementary *bars* with cross-sections $dx'dy'$ and boundaries $z'_{\min} = z'_{\min}(x', y')$ and $z'_{\max} = z'_{\max}(x', y')$ (roof and bottom according to the terminology of [31]) (figure 3).

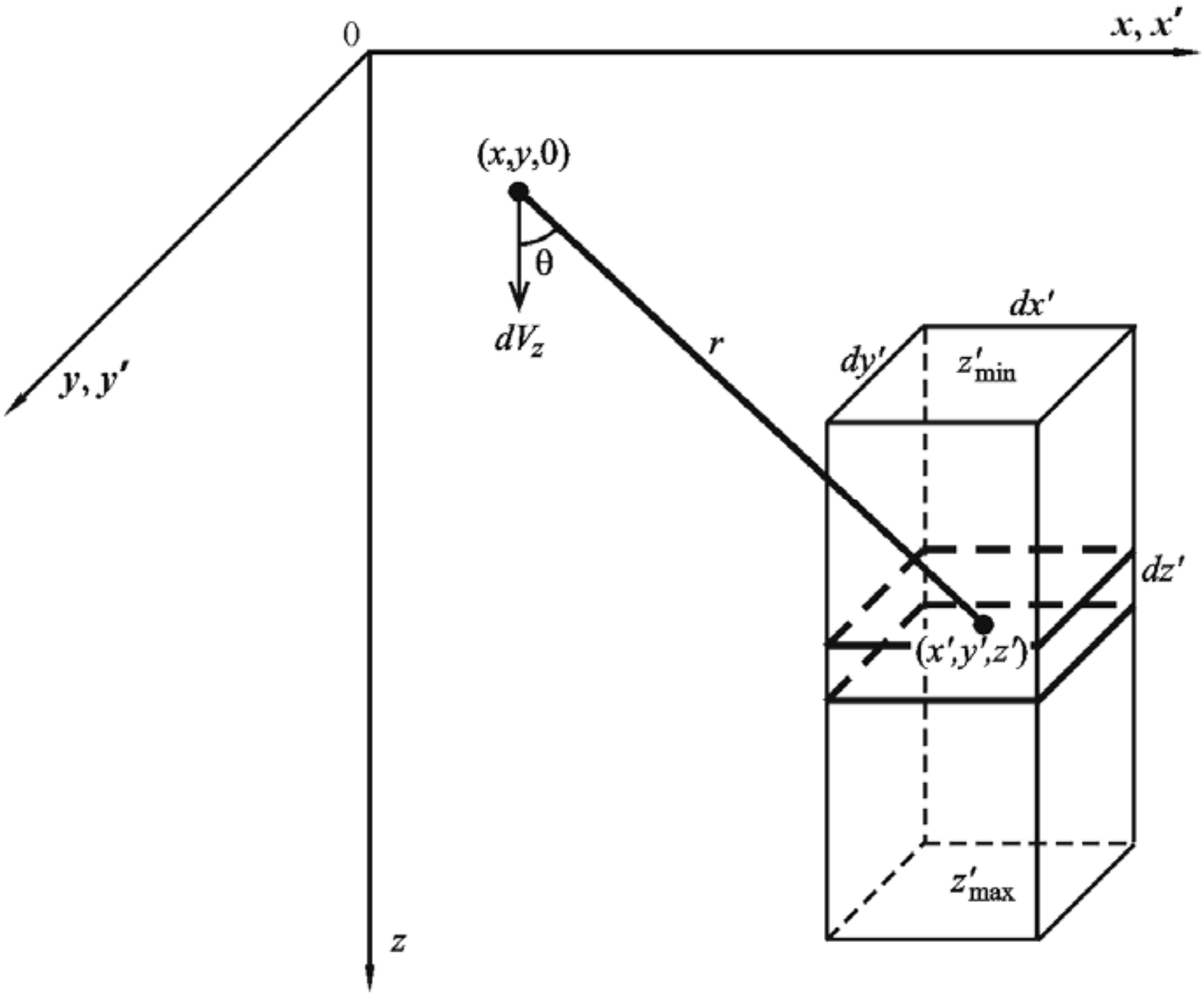

**Figure 3.** Elementary cell $dx'dy'dz'$ and elementary vertical bar of a body.

In the direct (and inverse) problem, we will consider only *z*-components of gravitational intensity of bodies. The intensity $V_z$ produced by an elementary bar at the point $(x, y, 0)$ is [8]

$$dV_z(x, y, 0) = \gamma\rho \left[ \frac{1}{\sqrt{(x-x')^2 + (y-y')^2 + z'^2_{\min}}} - \frac{1}{\sqrt{(x-x')^2 + (y-y')^2 + z'^2_{\max}}} \right] dx'dy', \quad (2.1)$$

where $\gamma$ is the gravitational constant, $\rho$ is the body density, and $x', y', z'$ are coordinates of a cell.

In order to calculate the intensity $V_z(x, y, 0)$ generated by the entire body at some point $(x, y, 0)$, it is necessary to integrate (2.1) over all elementary bars adjacent to each other from $z'_{\min} = z'_{\min}(x', y')$ to $z'_{\max} = z'_{\max}(x', y')$. If the condition for vertical star-shapedness is violated, one must exclude 'voids' from the integration regions $[z'_{\min}, z'_{\max}]$.

**Remark.** Note that in the direct problem, we do not approximate the body by a collection of vertical bars, but use them only to calculate the field. Furthermore, the body can have an enough arbitrary shape.

This technique for modeling the direct problem is sufficiently simple and effective, which is confirmed by solving numerical examples (see examples in [7, 8] and in this paper later).

## 3. Modeling the inverse problem with use of spheroids

**Definition 3.1.** An *ellipsoid* is a body bounded by the surface

$$\xi^2/a^2 + \eta^2/b^2 + \zeta^2/c^2 = 1,$$

where *a*, *b*, and *c* are semiaxes of the ellipsoid and the origin of the coordinate system is placed at the center of the body.

We will consider a biaxial ellipsoid, or ellipsoid of revolution about axis *z*. For brevity, we will call it a *spheroid*. Formulas for the potential *V* and the field components $V_x$, $V_y$, $V_z$ of a spheroid are deduced in [19–21]. However, these formulas are not brought to a final convenient form in the indicated works that we intend to do further for $V_z$.

**Remark.** As is commonly [19, 20], we will call an ellipsoid (spheroid) both the body bounded by the surface and the surface itself.

Consider an *oblate spheroid*, for which $a = b > c$. In [19–21], the coordinate system *x*, *y*, *z* is introduced with the origin at the center of the spheroid and with the *z*-axis directed vertically upwards along the minor axis *c* of the spheroid. For this case, formulas for the potential $V(x, y, z)$ and the intensity components $V_x(x, y, z)$, $V_y(x, y, z)$, $V_z(x, y, z)$ of a spheroid at some exterior point $(x, y, z)$ are deduced.

However, we consider the case, when the origin of the coordinate system *x*, *y*, *z* is placed at a point of the Earth's daily surface and the *z*-axis is directed vertically downward. Let the spheroid's center has the coordinates $x_0, y_0, z_0$ and the field (only the intensity $V_z$) of the spheroid be measured at the point $(x, y, 0)$ (figures 3 and 4).

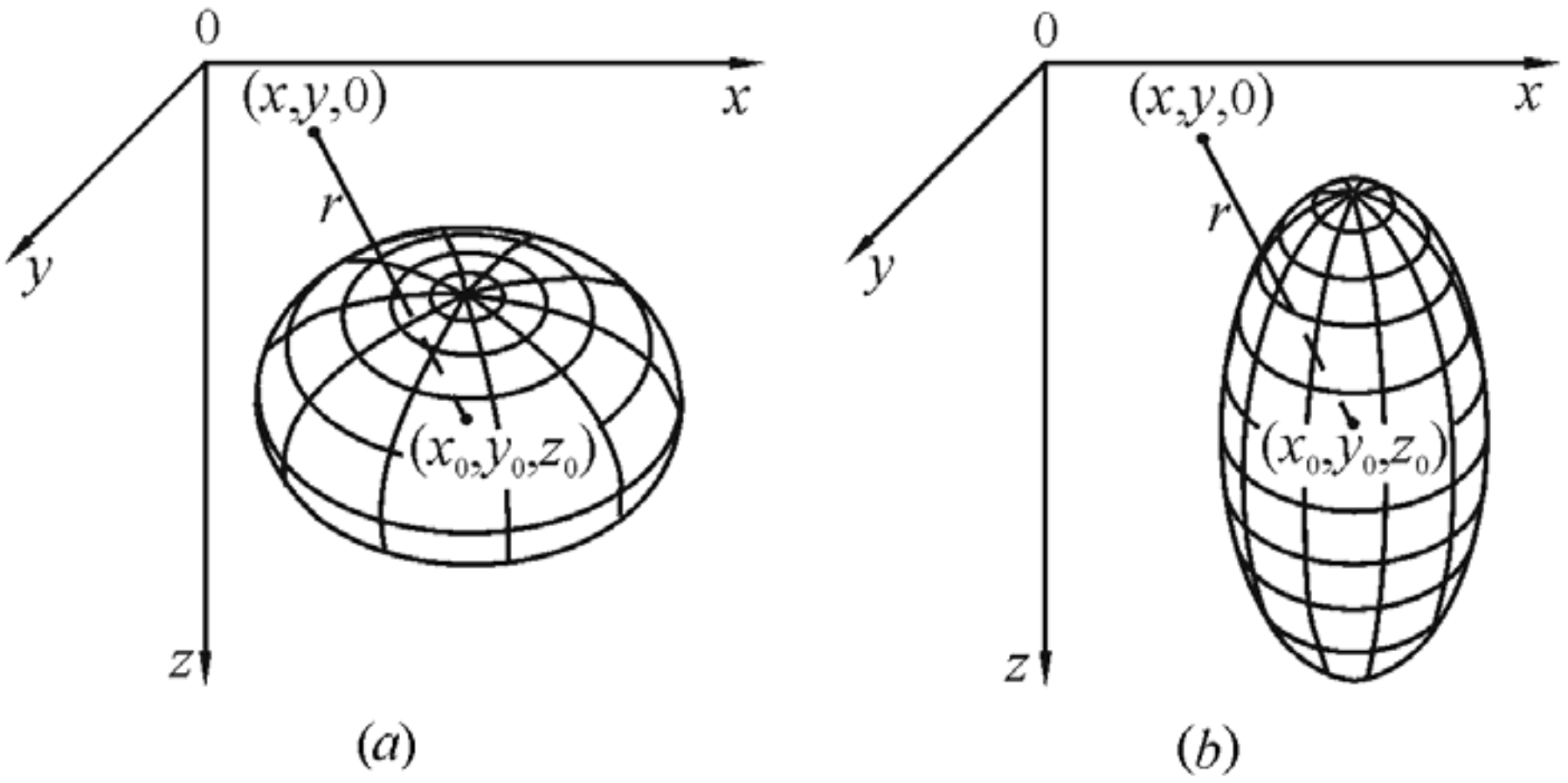

**Figure 4.** Oblate (*a*) and prolate (*b*) spheroids.

### (a) Formulas for *z*-intensity

For *oblate* (along *z*) *spheroid* (figure 4*a*), the formula for *z*-intensity takes the form (outside of spheroid) [8] (cf. [19–21]):

$$V_z(x, y, 0) = 4\pi\gamma\rho\frac{\varepsilon}{e^3}(p - \operatorname{arctg} p)\, z_0, \tag{3.1}$$

where

$$\varepsilon = c/a < 1,\quad e = \sqrt{1-\varepsilon^2} > 0,\quad p = q/\sqrt{\tau},\quad q = ea/r,\quad r = \sqrt{(x-x_0)^2 + (y-y_0)^2 + z_0^{\,2}},$$

$$\tau = \left[1 - q^2 + \sqrt{(1-q^2)^2 + 4q^2 z_0^{\,2}/r^2}\,\right]\Big/2. \tag{3.2}$$

For *prolate* (along *z*) *spheroid* (figure 4*b*), formula takes the form [8] (cf. [19–21]):

$$V_z(x,y,0) = 4\pi\gamma\rho\frac{\varepsilon}{e^3}\left[\ln\left(p+\sqrt{1+p^2}\right)-\frac{p}{\sqrt{1+p^2}}\right]z_0, \tag{3.3}$$

where

$$\varepsilon = c/a > 1,\ \ e=\sqrt{\varepsilon^2-1}>0,\ \ p=q/\sqrt{t},\ \ q=ea/r,\ \ r=\sqrt{(x-x_0)^2+(y-y_0)^2+z_0{}^2},$$

$$t=\left\{1-q^2+\sqrt{(1-q^2)^2+4q^2[(x-x_0)^2+(y-y_0)^2]/r^2}\right\}/2. \tag{3.4}$$

In the case of a sphere ($\varepsilon = c/a = 1$), the formula for *z*-intensity takes the form:

$$V_z(x,y,0) = \frac{4}{3}\pi\gamma\rho a^3\frac{z_0}{r^3}. \tag{3.5}$$

In this paper, we use only formulas (3.1)–(3.5) for *z*-intensity $V_z(x,y,0)$ of the oblate and prolate spheroids and the sphere.

The use of spheroids for solving the inverse gravimetry problem makes this approach close to the Sretenskii approach. Recall that according to the Sretenskii class models, a body possesses the mean plane *P* if any straight line perpendicular to this plane intersects the body surface only in two points on different sides of the plane *P*. If a geological body possesses a mean plane, its gravity center is inside the body, and its density is constant (and is given!), then the inverse problem (determining the body shape from the potential) has a unique solution. Our approach extends the Sretenskii approach because (see section 5) we do not assume that the density ρ is precisely given, but include it into the desired parameters (along with parameters $a$, ε, $x_0$, $y_0$, $z_0$). Furthermore, the solution uniqueness ensures at the expense of introducing constraints on the parameters (see later).

## 4. Initial approximations for parameters of a deposit

Let there be $m \geq 1$ bodies. Figure 5 plots isolines of the *z*-intensity $V_z(x,y,0)$ produced by bodies on the Earth's surface. Let $V_z(x_i,y_i,0)$ be measured at *N* points $i = 1,\ldots,N$. The number of bodies *m* and their coordinates $(x_0,y_0)_1,\ldots,(x_0,y_0)_m$ can be determined from the contour pattern (isolines).

### (a) Selection of bodies from isolines

The following w a y of a selection of bodies from the contour pattern (isolines) is proposed.

According to this way, *two conditions* must be fulfilled.

1. Valley in intensity $V_z$ between some two maxima (poles, peaks, hills) is not less than $\approx 20\%$ of intensities in poles (as in the Rayleigh criterion [27]).

2. The noise level $\delta V_z$ does not exceed 20% of intensities in poles[1].

If both conditions are fulfilled, we assume that two peaks (and hence two bodies) are determined from isolines. For example, figure 5 shows two poles with intensities $V_{z1} \approx 22$ and $V_{z2} \approx 30$. We assume $V_z = (V_{z1}+V_{z2})/2 \approx 26$. The intensity between the poles is $\upsilon_z \approx 17$, i.e. the valley is $(V_z - \upsilon_z)/V_z = 0.346 \approx 35\%$. In this case, the noise level is $\delta V_z \approx 5\%$ (see section 6). As a result, both conditions are fulfilled and we can assume that two bodies are determined from isolines (see also figure 8 later).

[1] The value 20% may vary. In addition, if the intensities in two poles are not equal: $V_{z1} \neq V_{z2}$, then we will use the value $V_z = (V_{z1}+V_{z2})/2$.

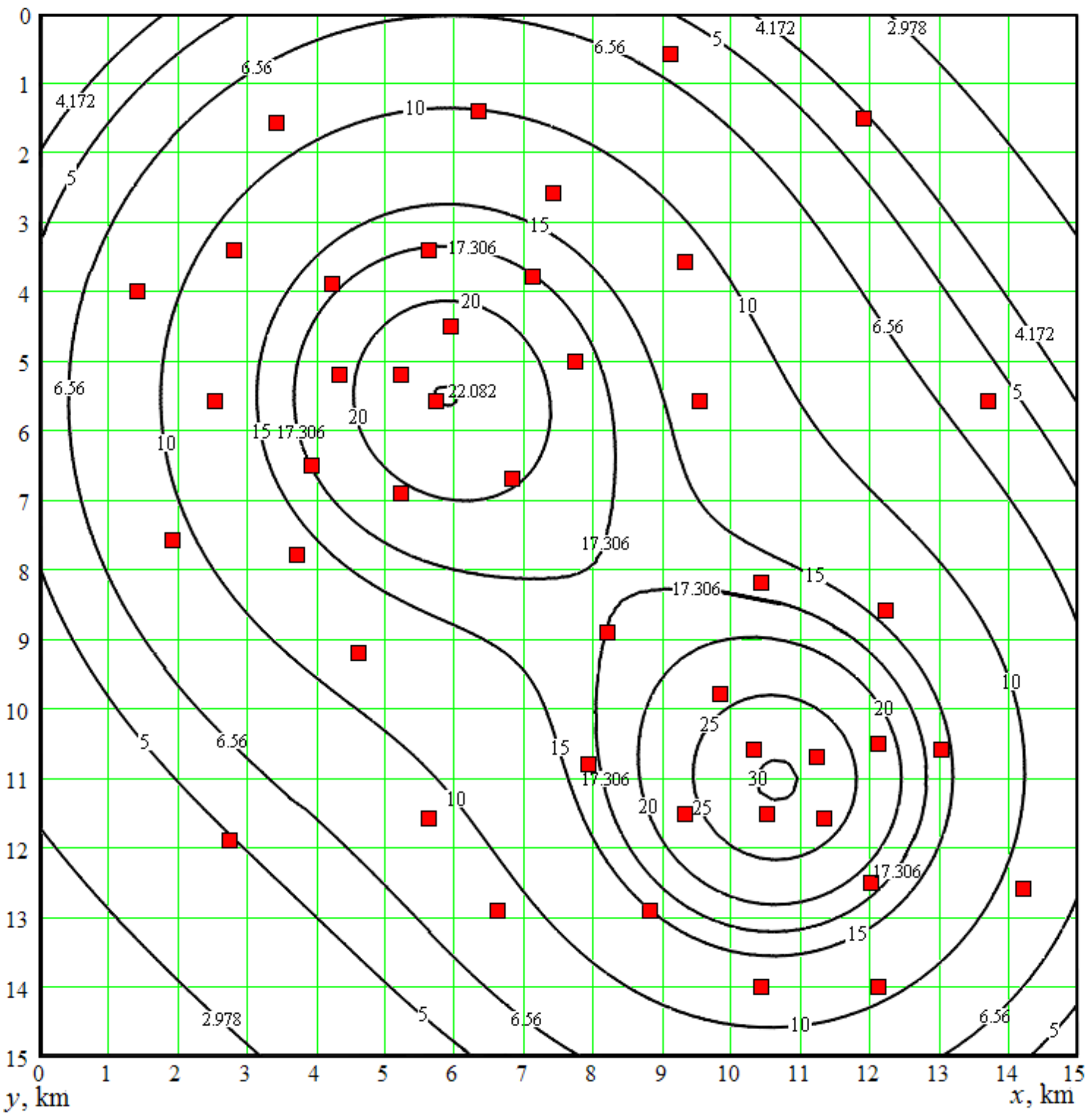

**Figure 5.** Isolines of intensity $V_z(x, y, 0)$, mGal from two bodies. (Online version in colour.)

The valley of 20% corresponds to some minimum distance between bodies in which they are delimited without mathematical processing. If the distance between the poles in figures 5 and 8 was 4 km, the valley would be equal to 20%, i.e. 4 km is the maximum distance between the poles in which they are resolved. If the valley is less than 20%, the bodies can be resolved mainly mathematically.

After separation of the bodies, we make an estimate of some parameters of the deposit.

### (b) Generalization of the Bulakh algorithm for estimating the depth and mass of a deposit

In [32], an algorithm for estimating the depth $z_0$ and mass $M$ of each body was proposed. We present this algorithm (the *Bulakh algorithm*). We assume for a time that there is only one body, for example, a body matched by the $V_z$ isolines in the upper left part of figure 5. We assume also that this body is a *homogeneous sphere*. Denote by $R(x_0, y_0, z_0)$ its center, $M$ its mass, and $Q(x_0, y_0, 0)$ the point at which $V_z = \max$ (figure 6). Then, the $z$-intensity $V_z(x, y, 0)$ at some point $P(x, y, 0)$ is

$$V_{zP} \equiv V_z(x, y, 0) = V_z(s) = \gamma \frac{M}{r^2} \cdot \frac{z_0}{r} = \gamma \frac{M z_0}{(z_0^2 + s^2)^{3/2}}, \tag{4.1}$$

where $r = \sqrt{z_0^2 + s^2}$ is the distance between $R$ and $P$ and $s = \sqrt{(x - x_0)^2 + (y - y_0)^2}$ is the distance between $Q$ and $P$. The $z$-intensity at the point $Q$ is

$$V_{zQ} \equiv V_z(x_0, y_0, 0) = V_{z\max} = \gamma \frac{M}{z_0^2}. \tag{4.2}$$

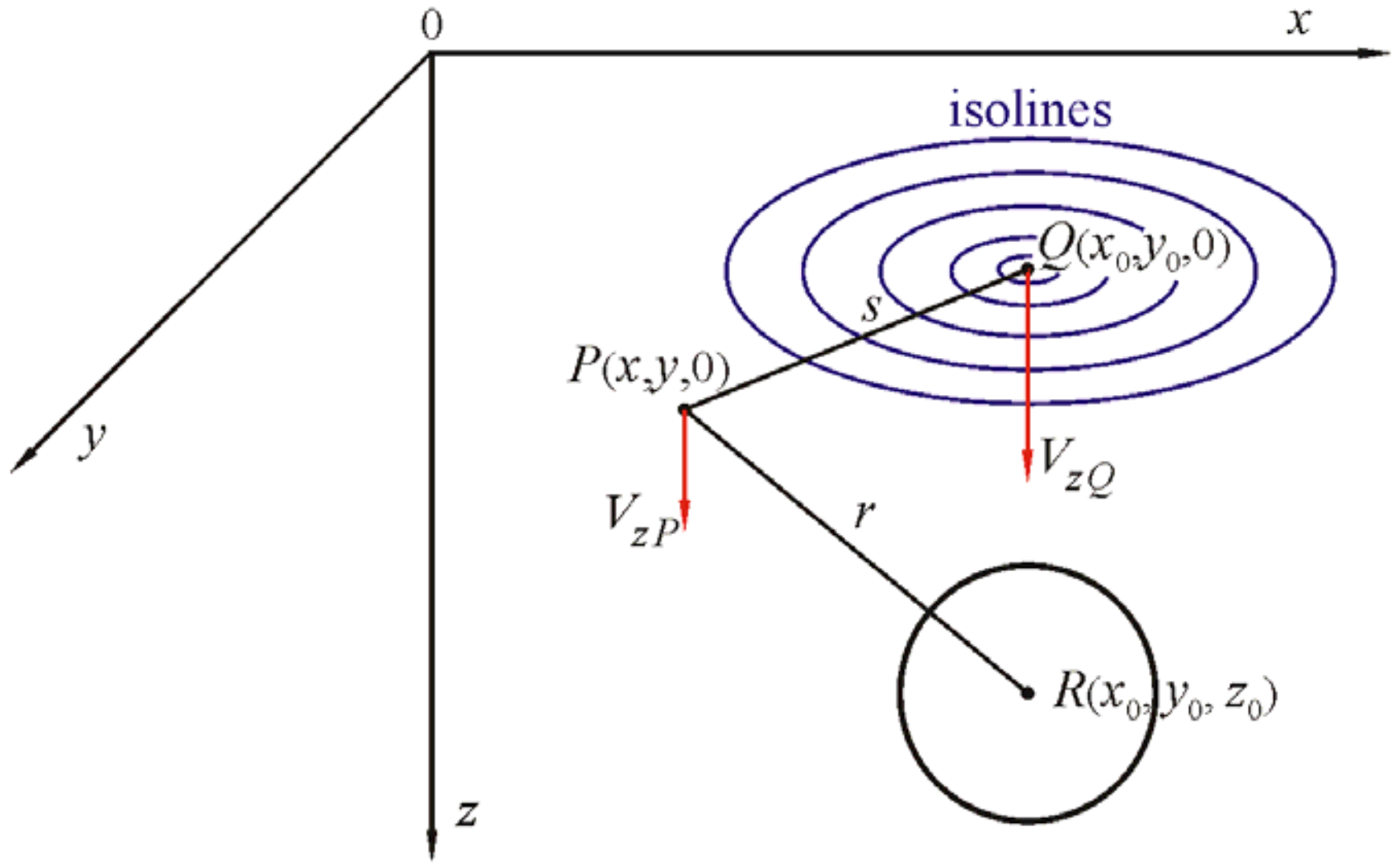

**Figure 6.** Deposit in the form of one body (homogeneous sphere). (Online version in colour.)

Relations (4.1) and (4.2) can be considered as a system of two equations for $z_0$ and *M*. Using the notation $\nu = V_{zP}/V_{zQ}$, we obtain $\nu = (z_0/s)^3/[(z_0/s)^2+1]^{3/2}$ or

$$\left(\frac{\mu^2}{\mu^2+1}\right)^{3/2} = \nu, \tag{4.3}$$

where $\mu = z_0/s$. Relation (4.3) is equation for $\mu = \mu(\nu)$ at given (measured) ν. Its solution is

$$\mu(\nu) = \sqrt{\frac{\nu^{2/3}}{1-\nu^{2/3}}}, \quad \nu \in [0,1], \quad \mu \in (0,\infty). \tag{4.4}$$

In table 1, the dependence $\mu(\nu)$ is presented.

**Table 1.** Dependence $\mu(\nu)$

| $\nu = \frac{V_z(s)}{V_{z\max}}$ | 0 | 0.1 | 0.2 | 0.3 | 0.4 | 0.5 | 0.6 | 0.7 | 0.8 | 0.9 | 1 |
|---|---|---|---|---|---|---|---|---|---|---|---|
| $\mu = \frac{z_0}{s}$ | 0 | 0.5240 | 0.7209 | 0.9011 | 1.0898 | 1.3048 | 1.5700 | 1.9301 | 2.4969 | 3.7071 | ∞ |

Using this dependence, we can estimate the depth $z_0 = \mu s$ of a body from measured ν and *s*. It is desirable that the point *P* is located away from other bodies, e.g., at $x < 7$ km and $y < 6$ km in figure 5. One can make several estimates of μ (and $z_0$) from a few points *P* and average the result. Then this procedure must be performed for each of the *m* bodies.

**Remark.** Definition of the depth $z_0$ and mass *M* of a body is possible only if $s \neq 0$, i.e. when measurements are performed at two different points *P* and *Q*. Moreover, to improve the accuracy of the algorithm it is necessary that $s \approx z_0$, as an analysis of formulas (4.1)–(4.4) and solving numerical examples shows.

This algorithm is applicable only if a measurement is performed in the maximum value of $V_z$ as in figure 5 for the upper left body. If such a measurement is not performed as in the lower right corner of figure 5, but $x_0, y_0$ can be estimated, then we make a *generalization of the Bulakh algorithm*. As a point $Q(x_Q, y_Q, 0)$ we take a point with value of $V_z$, close to $V_{z\max}$. Then

$$V_{zQ} \equiv V_z(x_Q, y_Q, 0) = \gamma \frac{M}{d^2} \cdot \frac{z_0}{d} = \gamma \frac{M z_0}{\left(z_0^2 + \Delta^2\right)^{3/2}}, \tag{4.5}$$

where $d = \sqrt{z_0^2 + \Delta^2}$ is the distance between *Q* and *R*, $\Delta = \sqrt{(x_Q - x_0)^2 + (y_Q - y_0)^2}$ is the distance between *Q* and the point $(x_0, y_0, 0)$, and expression (4.1) remains valid for $V_{zP}$.

Formula (4.5) generalizes formula (4.1). The ratio of the intensities at points *P* and *Q* is

$$\nu = \frac{V_{zP}}{V_{zQ}} = \left(\frac{\mu^2 + \psi^2}{\mu^2 + 1}\right)^{3/2},$$

where $\psi = \Delta/s$, or

$$\left(\frac{\mu^2(\nu,\psi) + \psi^2}{\mu^2(\nu,\psi) + 1}\right)^{3/2} = \nu. \tag{4.6}$$

Relation (4.6) is the equation over $\mu = \mu(\nu,\psi)$ at given $\nu$ and $\psi$. It generalizes equation (4.3). Its solution is

$$\mu(\nu,\psi) = \sqrt{\frac{\nu^{2/3} - \psi^2}{1 - \nu^{2/3}}}, \quad \nu \in [0,1], \quad \psi \in (0,1), \quad \mu \in (0,\infty). \tag{4.7}$$

In figure 7, the dependence $\mu(\nu,\psi)$ is ploted.

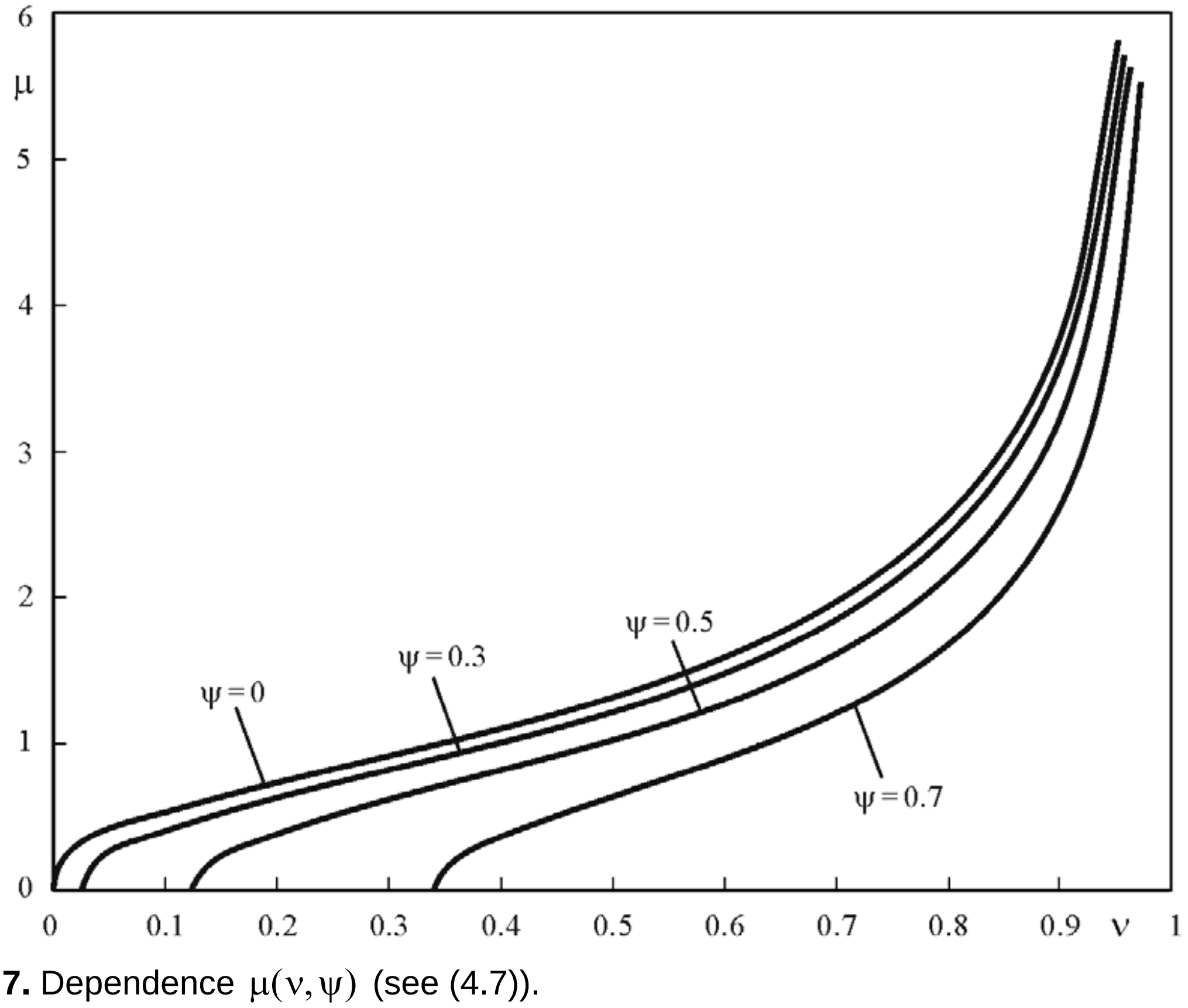

**Figure 7.** Dependence $\mu(\nu,\psi)$ (see (4.7)).

After estimating μ by (4.7), one can estimate the *depth* of the mass center $z_0 = \mu s$ for the body. Such estimate should be made for each of the *m* bodies. Now we can estimate the *mass* *M* of the body. If there is a measurement of $V_z$ at the point $(x_0, y_0, 0)$ of its maximum, we obtain from (4.2)

$$M = \frac{1}{\gamma} z_0^2 V_{z\,\max}\,. \tag{4.8}$$

If there is no $V_z$ at the point $(x_0, y_0, 0)$ of its maximum, then using a measurement of $V_z$ at some point $Q(x_Q, y_Q, 0)$ we obtain from (4.5):

$$M = \frac{1}{\gamma} \frac{(z_0^2 + \Delta^2)^{3/2}}{z_0} V_{zQ}\,. \tag{4.9}$$

At $\Delta = 0$, formula (4.9) turns into (4.8).

Let the coordinates and distances *x*, *y*, *z*, *s*, *r*, Δ, *d* are expressed in kilometers, the mass *M* in billions of tons, and a gravitational intensity (anomaly) $V_z$ in mGal. Then formulas (4.8) and (4.9) take the form

$$M = 0.15 z_0^2 V_{z\,\max}\,, \tag{4.10}$$

$$M = 0.15 \frac{(z_0^2 + \Delta^2)^{3/2}}{z_0} V_{zQ}\,. \tag{4.11}$$

These results generalize the algorithm adduced in [32]. They give good initial approximations to $z_0$ and *M* (as well as $x_0$ and $y_0$) for each body. In the next section, we show how to make more precise these (and other) parameters of nonspherical (spheroidal) bodies.

## 5. Refinement of the deposit parameters

Let us consider the problem for more precise determining the parameters of spheroids modeling bodies of a deposit. Let $\tilde{V}_{zi} \equiv \tilde{V}_z(x_i, y_i, 0)$, $i = 1, \ldots, N$ be values of $V_z$ measured (with errors) at a series of points $(x_i, y_i, 0)$, $i = 1, \ldots, N$ on the Earth's surface (in particular, on a ship), where *N* is the number of measurement points. Let $V_{zi} \equiv V_z(x_i, y_i, 0)$ denote values of $V_z$ calculated from formulas (3.1), (3.3) or (3.5), *k* denote the number of sought parameters of every from *m* spheroids (the types of parameters are listed later), in all *mk* parameters, and $p_j$, $j = 1, \ldots, mk$ denote the deposit parameters.

The problem for determining the deposit parameters is ill-posed, viz., unstable and having, generally speaking, a nonunique solution [1, 6–8, 30–34]. Even if a body is star-shaped and has the Sretenskii mean plane, a solution will be nonunique if the body density ρ is not given. Most effectively, an elimination of the solution nonuniqueness can be realized via using solution constraints and the Tikhonov regularization method.

We will solve this problem via *minimizing the Tikhonov smoothing functional* [16] (cf. [23, 28, 35]) in two variants [6–8]:

$$F_1 \equiv \sum_{i=1}^{N} \left(\tilde{V}_{zi} - V_{zi}\right)^2 + \alpha \sum_{j=1}^{mk} q_j (p_j - p_{\mathrm{mid}\,j})^2 = \min_{p_1, \ldots, p_{mk}} \tag{5.1}$$

or

$$F_2 \equiv \sum_{i=1}^{N} \left(\tilde{V}_{zi} - V_{zi}\right)^2 + \alpha \sum_{j=1}^{mk} q_j p_j^2 = \min_{p_1, \ldots, p_{mk}}\,, \tag{5.2}$$

where $q_j$ are weights and $\alpha > 0$ is the regularization parameter. In order to enhance the stability and eliminate the solution nonuniqueness, we impose *constraints on solution* in the form of prescribed ranges of parameter variations (cf. [34, 36]):

$$p_{\min j} \le p_j \le p_{\max j}, \quad j = 1, \ldots, mk, \tag{5.3}$$

where $p_{\min j}$ and $p_{\max j}$ must be estimated from an additional (a priori) information. We assume for the weights that $q_j = 1 \big/ p_{\mathrm{mid}\, j}^2$, where $p_{\mathrm{mid}\, j} = (p_{\min j} + p_{\max j})/2$. Such problem is a problem of nonlinear programming with regularization [16, 37–39].

**Remark.** We introduce the weights $q_j$ because the sought parameters can have a different physical dimension and a different order of magnitudes, and inroduction of the weights makes the summands $q_j (p_j - p_{\mathrm{mid}\, j})^2$ in (5.1) or $q_j p_j^2$ in (5.2) by close (and dimensionless) quantities.

The mean square error of the solution is [6, 8]

$$\delta = \left[ \frac{1}{mk} \sum_{j=1}^{mk} q_j (p_j - \bar{p}_j)^2 \right]^{1/2}, \tag{5.4}$$

where $p_j = p_j(\alpha)$, $\delta = \delta(\alpha)$, and $\bar{p}_j$ are 'exact' $p_j$ estimated from the modeled body.

We offer the following s c h e m e for calculating the parameters of deposit bodies.

1. The number of bodies *m* and values of $x_0$, $y_0$ for each body are estimated from the pattern of $V_z$ isolines (see figure 5).

2. Using measurements of $V_z$ at a series of points, we determine the averaged value of $\mu$ by formula (4.4) or/and (4.7) and estimate $z_0 = \mu s$ for each body.

3. The mass *M* of each body is estimated from formulas (4.8)–(4.11).

4. The parameters $p_1 = \varepsilon$, $p_2 = \rho$, $p_3 = x_0$, $p_4 = y_0$, and $p_5 = z_0$ ($k = 5$) are sought for each body, and the obtained estimations of $x_0, y_0, z_0$ are used as initial approximations by setting $p_{\mathrm{mid}\, 3} = x_0$, $p_{\mathrm{mid}\, 4} = y_0$, $p_{\mathrm{mid}\, 5} = z_0$. Furthermore, the estimate of *M* is used for calculating the semiaxis *a*,

$$a = \sqrt[3]{\frac{M}{(4/3)\pi\varepsilon\rho}}. \tag{5.5}$$

There exists a number of methods for minimizing functionals with constraints (of conditional gradient, conjugate gradients, steepest descent, ravines, chords, et al.) [16, 37, 38]. In exploration geophysics, the *method of gradient descent* is often used (works of Kantorovich, et al.) [31]. However, a point of convergence of this method strongly depends on initial approximation of a solution. Besides, this method uses derivatives of a functional and does not use constraints on a solution.

For minimizing the functional (5.1) or (5.2), in this paper we use the modification of the *coordinate descent method* [40] effective in the case of constraints of the form (5.3). Specific character of this method consists in the fact that it does not allow to the solution to fall outside the constraints (5.3), thus providing stability and convergence of the solution within of the interval given by the inequalities (5.3). However, it is not always possible to set a sufficiently narrow internal. In this case, one should set initial wide constraints and then set them iteratively narrower so that the solution does not go out from the interval (5.3). Just so, the following examples were solved. We call this method a *way of decremental constraints*.

## 6. Numerical examples (modeling the direct and inverse problems)

The *software packages* IGP2 and IGP5 (Inverse Gravimetry Problem, var. 2 and 5) were developed for modeling the direct and inverse gravimetry problems where a deposit consists of two or five bodies ($m = 2$ or $m = 5$). Calculations were performed using the programming language MS Fortran PowerStation 4.0 and graphs were constructed in MathCAD, CorelDRAW, and Paint. The following two *numerical examples* were solved.

### (a) Example 1

The deposit is prescribed in the form of two bodies (figures 1 and 2). In solving the direct problem (modeling the intensities $V_z$ created by the deposit), one sets $N = 45$ points $(x_i, y_i, 0)$, $i = 1,\ldots,N$ for measurement of $V_z$ on the Earth's surface (the small squares in figure 5). At each point $(x_i, y_i, 0)$, $i = 1,\ldots,N$, the intensity $V_z(x_i, y_i, 0)$ is modeled via summing over elementary bars and over two bodies using formula (2.1) with the discretization step $dx' = dy' = 0.25\,\text{km}$ for the first body and $dx' = dy' = 0.125\,\text{km}$ for the second body with prescribed $z'_{\min}(x', y')$ and $z'_{\max}(x', y')$ (see figure 3).

With the help of a random number generator [25], we add random errors $\delta V_z$ distributed normally with zero expectation and standard deviation $\sigma = 1$ mGal ($\approx 6\%$ of the average value of $V_z$ and $\approx 4\%$ of $V_{z\max}$) to exact values $V_z(x_i, y_i, 0)$. The errors $\delta V_z$ model measurement errors, small-scale nonhomogeneities of a medium, etc.

On the basis of values $V_z(x_i, y_i, 0)$ at $x \in [0,15]$ km, $y \in [0,15]$ km, isolines were constructed with the help of the comand Contour Plot in MathCAD. In figure 5, the isolines of intensity $V_z(x, y, 0)$ (mGal) are represented. Furthermore, $\rho = 1.6\,\text{g/cm}^3$ for the first body and $\rho = 2.6\,\text{g/cm}^3$ for the second body.

However in practice, continuous data $V_z(x, y, 0)$ are usually unknown and only (noisy) values of $\tilde{V}_z$ at measurement points $(x_i, y_i, 0)$, $i = 1,\ldots,N$, are known. In figure 8, isolines of z-intensity $\tilde{V}_z = V_z + \delta V_z$ (mGal) constructed from measured values $\tilde{V}_z(x_i, y_i, 0)$, $i = 1,\ldots,N$, are represented.

**Remark.** Strictly speaking, isolines of $V_z(x, y, 0)$ are not needed for constructing a deposit model. As is obvious from (5.1) and (5.2), only discrete values of $\tilde{V}_{zi}$ are used. Nevertheless, isolines are needed for estimating quantity of bodies and their coordinates $x_0$ and $y_0$.

Now we present the results of the solving the *inverse problem*, viz. determining the parameters for two bodies of the deposit.

First, we make estimates of $x_0, y_0$ for both bodies from the isolines of figure 8. For the first body, we obtain $x_0 = 5.7$ and $y_0 = 5.6$ km. Estimates of $z_0 = \mu s$ using formulas (4.4) and (4.7) by a few points $P$ and $Q$ in the area $x \le 7$ & $y \le 6$ km are obtained from 3 to 5.5 km; on average, $z_0 = 4.9$ km. Estimates of mass $M$ by formulas (4.8) and (4.9) (or (4.10) and (4.11)) are obtained from 30 to 105 bln t; on average, $M = 67.5$ bln t.

For the second body, $x_0 = 10.65$ and $y_0 = 11.1$ km. Estimates of $z_0 = \mu s$ using formula (4.7) by a few points for $x \ge 8$ & $y \ge 10$ km are obtained from 2.9 to 3.7 km; on average, $z_0 = 3.3$ km. Estimates of mass $M$ by formula (4.11) are obtained from 37 to 60 bln t; on average, $M = 48.5$ bln t.

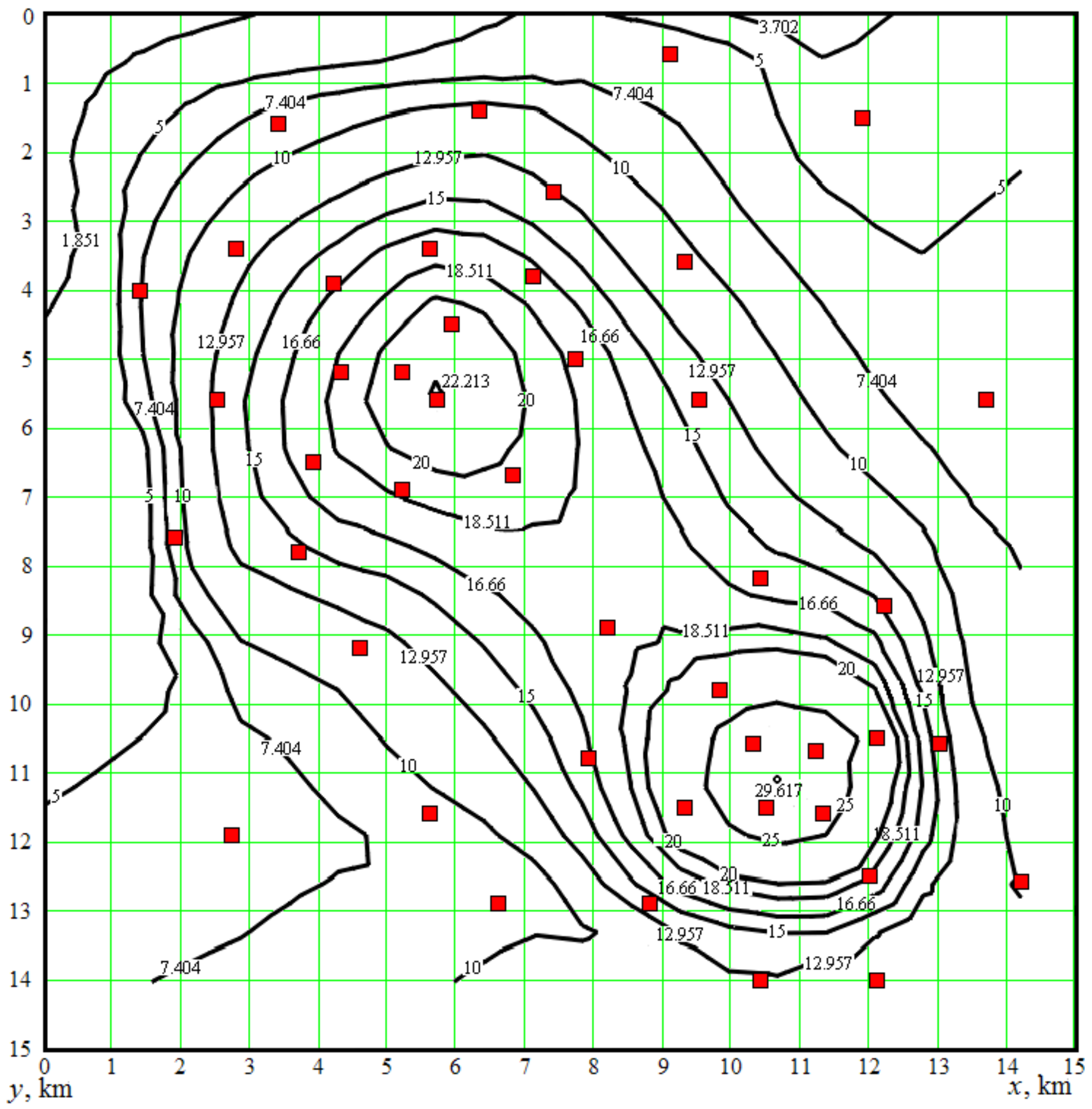

**Figure 8.** Isolines of intensity $\widetilde{V}_z(x, y, 0)$, mGal from two bodies constructed by $N = 45$ points of measurements (example 1). (Online version in colour.)

Then we determine (or make more precise) the following ten sought parameters via minimizing the functional $F_1$ or $F_2$ (see (5.1) or (5.2)): $p_1 = \varepsilon$, $p_2 = \rho$, $p_3 = x_0$, $p_4 = y_0$, $p_5 = z_0$ for the first body and $p_6 = \varepsilon$, $p_7 = \rho$, $p_8 = x_0$, $p_9 = y_0$, $p_{10} = z_0$ for the second body. In table 2, the initial lower constraints $p_{\min}$ and upper constraints $p_{\max}$ and the arithmetic mean values $p_{\text{mid}} = (p_{\min} + p_{\max})/2$ are presented. Such (approximate) values can be given by practiced specialists and then they are refined (the constraints are made narrower) so that the sought solution $p_i, \ldots, p_{10}$ does not go beyond the interval (5.3). Furthermore, estimate of the semiaxis $a$ is realized by formula (5.5).

The dependence of the mean square error $\delta$ of the solution (see (5.4)) on the regularization parameter $\alpha$ is calculated. Furthermore, the following parameters of the bodies-spheroids based on the deposit estimation (figures 1 and 2) are taken as 'exact': $\bar{p}_1 = \bar{\varepsilon} = 0.51$, $\bar{p}_2 = \bar{\rho} = 1.6\ \text{g/cm}^3$, $\bar{p}_3 = \bar{x}_0 = 5.7\ \text{km}$, $\bar{p}_4 = \bar{y}_0 = 5.3\ \text{km}$, $\bar{p}_5 = \bar{z}_0 = 4.2\ \text{km}$ for the first body and $\bar{p}_6 = \bar{\varepsilon} = 1.96$, $\bar{p}_7 = \bar{\rho} = 2.6\ \text{g/cm}^3$, $\bar{p}_8 = \bar{x}_0 = 10.7\ \text{km}$, $\bar{p}_9 = \bar{y}_0 = 11.1\ \text{km}$, $\bar{p}_{10} = \bar{z}_0 = 3.8\ \text{km}$ for the second body.

Since a regularization is used, so it is necessary to turn our attention to the question about choise of the regularization parameter $\alpha$. There exists a number of the ways for choosing $\alpha$ (the discrepancy principle, et al.) [16–18, 25]. However, in example 1, rigid constraints are imposed on the solution $p$ and they provide stability and uniqueness of the solution without regularization. The problem has practically the same solution for any $\alpha \le 10^{-8}$ and the solu-tion does not depend on the type of the functional ($F_1$ or $F_2$). But if constraints are less rigid than in table 2, the error $\delta(\alpha)$ and the solution $p(\alpha)$ depend on $\alpha$ and the error $\delta$ at some $\alpha = \alpha_{opt}$ may be smaller than at $\alpha = 0$.

**Table 2.** Refined constraints on parameters and obtained solution

| $p$ | $\varepsilon$ | $\rho$, $g/cm^3$ | $x_0$, km | $y_0$, km | $z_0$, km | $\upsilon$, $km^3$ | $M$, bln t |
|---|---|---|---|---|---|---|---|
| First body | | | | | | | |
| $p_{min}$ | 0.2 | 1.1 | 5.4 | 5.2 | 4.0 | | |
| $p_{max}$ | 0.6 | 1.7 | 6.0 | 6.0 | 5.8 | | |
| $p_{mid}$ | 0.4 | 1.4 | 5.7 | 5.6 | 4.9 | | |
| solution | 0.595 | 1.69 | 5.75 | 5.48 | 4.40 | 39.86 | 67.50 |
| Second body | | | | | | | |
| $p_{min}$ | 1.8 | 2.3 | 10.3 | 10.2 | 2.3 | | |
| $p_{max}$ | 2.2 | 2.9 | 11.0 | 12.0 | 4.3 | | |
| $p_{mid}$ | 2.0 | 2.6 | 10.65 | 11.1 | 3.3 | | |
| solution | 2.035 | 2.65 | 10.71 | 11.90 | 3.81 | 18.28 | 48.50 |

In table 2, the obtained solution is given, viz. the values of the parameters $\varepsilon$, $\rho$, $x_0$, $y_0$, $z_0$, $\upsilon$ and $M$ for both bodies-spheroids. For the first body, $a = 2.52$ km, $\varepsilon = 0.595$, $\rho = 1.69\ g/cm^3$, $x_0 = 5.75$ km, $y_0 = 5.48$ km, $z_0 = 4.40$ km. For the second body, $a = 1.29$ km, $\varepsilon = 2.035$, $\rho = 2.65\ g/cm^3$, $x_0 = 10.71$ km, $y_0 = 11.90$ km, $z_0 = 3.81$ km. We see that all parameters are obtained within the interval (5.3). Furthermore, the initial geologic deposit has the following parameters (see figures 1 and 2).

For the first body, the half-length along the *x* direction (*x* half-length) is 2.5 km and the *y* half-length is 3 km, i.e., the average half-length in the plane *xy* is 2.75 km; the *z* half-length is 1.4 km, i.e., the ratio between the *z* и *xy* half-lengthes is 0.51; the coordinates of the center are $x_0 = 5.7$ km, $y_0 = 5.3$ km, $z_0 = 4.2$ km; the density $\rho$ is $1.6\ g/cm^3$. For the second body, the *x* half-length is 1.25 km, the *y* half-length is 1.5 km, the average *xy* half-length is 1.375 km; the *z* half-length is 2.695 km, the ratio between the *z* и *xy* half-lengthes is 1.96; the coordinates of the center are $x_0 = 10.7$ km, $y_0 = 11.1$ km, $z_0 = 3.8$ km; the density $\rho$ is $2.6\ g/cm^3$.

In figure 9, the images of two calculated spheroids are given. We see that the modeling of deposit by spheroids gives a quite satisfactory result. Note that the initial value of the functional $F_1$ (before its minimization) is $F_1 = 0.167 \cdot 10^{-2}$, and after minimization by the coordinate descent method $F_1 = 0.197 \cdot 10^{-3}$, i.e. the value of $F_1$ decreases by one order.

Volumes of two spheroids $\upsilon = (4/3)\pi a^3 \varepsilon$ and their masses $M = \upsilon\rho$ are also calculated. For the first spheroid, $\upsilon = 39.86\ km^3$, $M = 67.50$ bln t, and for the second spheroid, $\upsilon = 18.28\ km^3$, $M = 48.50$ bln t. For comparison, the true volumes $\overline{\upsilon}$ and masses $\overline{M}$ of the deposit bodies (figure 1) calculated by summation over bars are $\overline{\upsilon} = 39.67\ km^3$, $\overline{M} = 63.48$ bln t for the first body and $\overline{\upsilon} = 19.06\ km^3$, $\overline{M} = 49.55$ bln t for the second body.

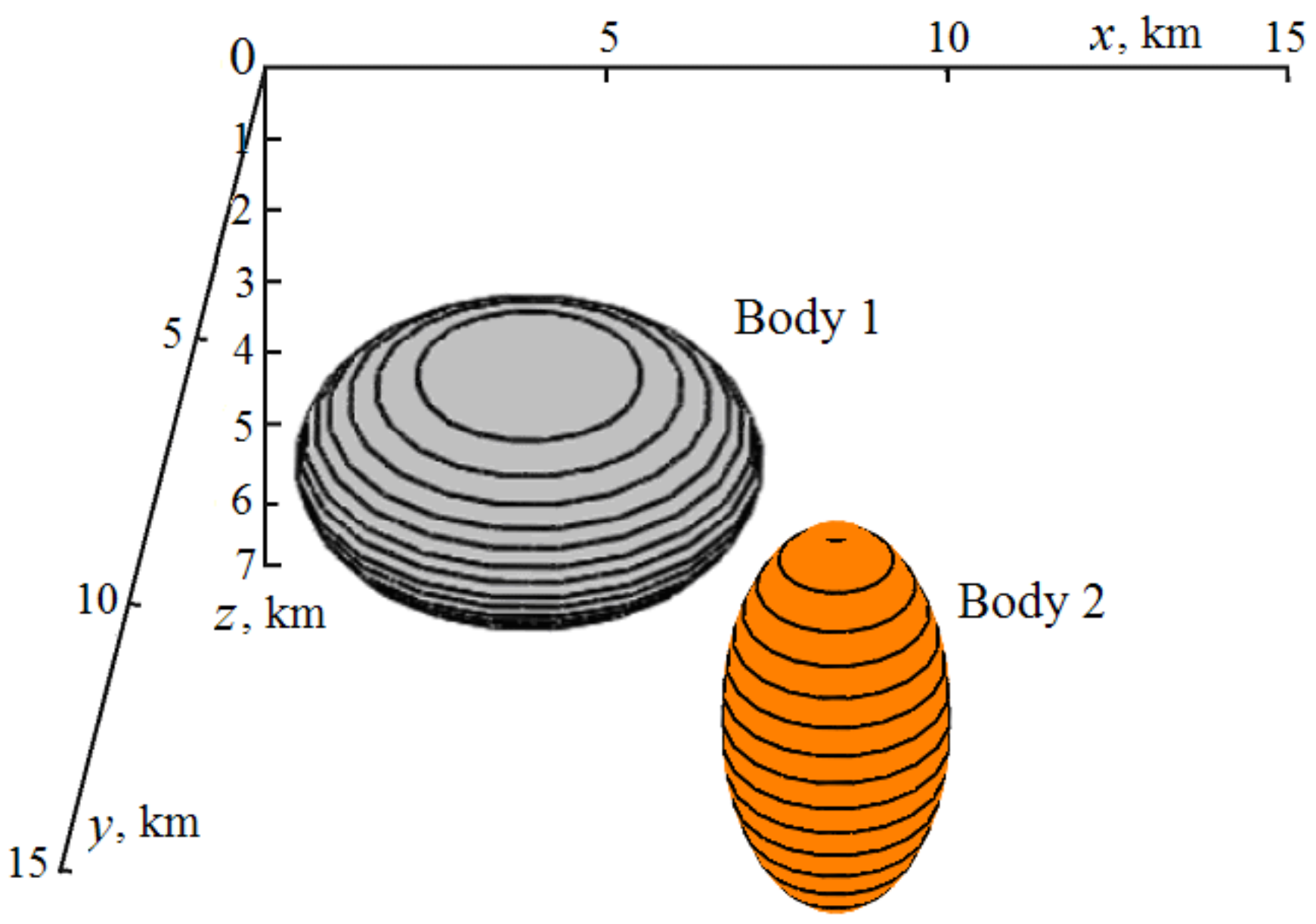

**Figure 9**. Two spheroids imitating the deposit in example 1. (Online version in colour.)

## (b) Example 2

This (more complicated) example is characterized by the following.

(1) The deposit consists of a larger number of bodies than in example 1.

(2) As the input data, noisy experimental values of $z$-intensities $\widetilde{V}_z(x_i, y_i, 0)$, $i = 1, \ldots, N$, in $N = 73$ points on the Earth's surface (squares in figure 10) are given.

(3) Exact values of parameters of bodies are unknown. Therefore, in order that to estimate the accuracy of determining the parameters of bodies-spheroids, we will use a comparison of isolines constructed by values of $\widetilde{V}_z$ (figure 10) and ones constructed by results of calculating the spheroids (see figure 11 later).

On the bases of values of $\widetilde{V}_z$ and with the help of computerized graphics of the system MathCAD, isolines of $z$-intensity are constructed (figure 10).

By the isolines of figure 10 using the way of a selection of bodies from the isolines (see section 4a), we conclude that the number of bodies in deposit is $m = 5$. By the isolines, it is possible to estimate the coordinates $(x_0, y_0)$ of every body. With the help of generalization of the Bulakh algorithm and with use of formulas (4.4), (4.7), (4.10), (4.11), we obtain also the estimations of the depths $z_0 = \mu s$ and masses $M$ of bodies. Results are the following ($x_0, y_0, z_0$ in km, $M$ in bln t):

for body 1: $x_0 = 3$, $y_0 = 3.5$, $z_0 = 4.45$, $M = 55.9$;

for body 2: $x_0 = 10.3$, $y_0 = 11.3$, $z_0 = 4.1$, $M = 47.5$;

for body 3: $x_0 = 3$, $y_0 = 11.8$, $z_0 = 3.8$, $M = 37.5$;

for body 4: $x_0 = 11$, $y_0 = 2$, $z_0 = 4.4$, $M = 43.1$;

for body 5: $x_0 = 12.9$, $y_0 = 6.3$, $z_0 = 5.1$, $M = 30$.

Further, via minimizing the functional $F_1$ (see (5.1)), one corrects (determines) the sought parameters of the bodies-spheroids $p_1 = \varepsilon$, $p_2 = \rho$, $p_3 = x_0$, $p_4 = y_0$, $p_5 = z_0$ for body 1 and analogously for bodies 2–5, in all $mk = 25$ parameters. Furthermore, at first wide con-straints on values of the parameters are prescribed and later they become by more narrow so that the sought solution $p_1, \ldots, p_{25}$ does not go out from the interval (5.3).

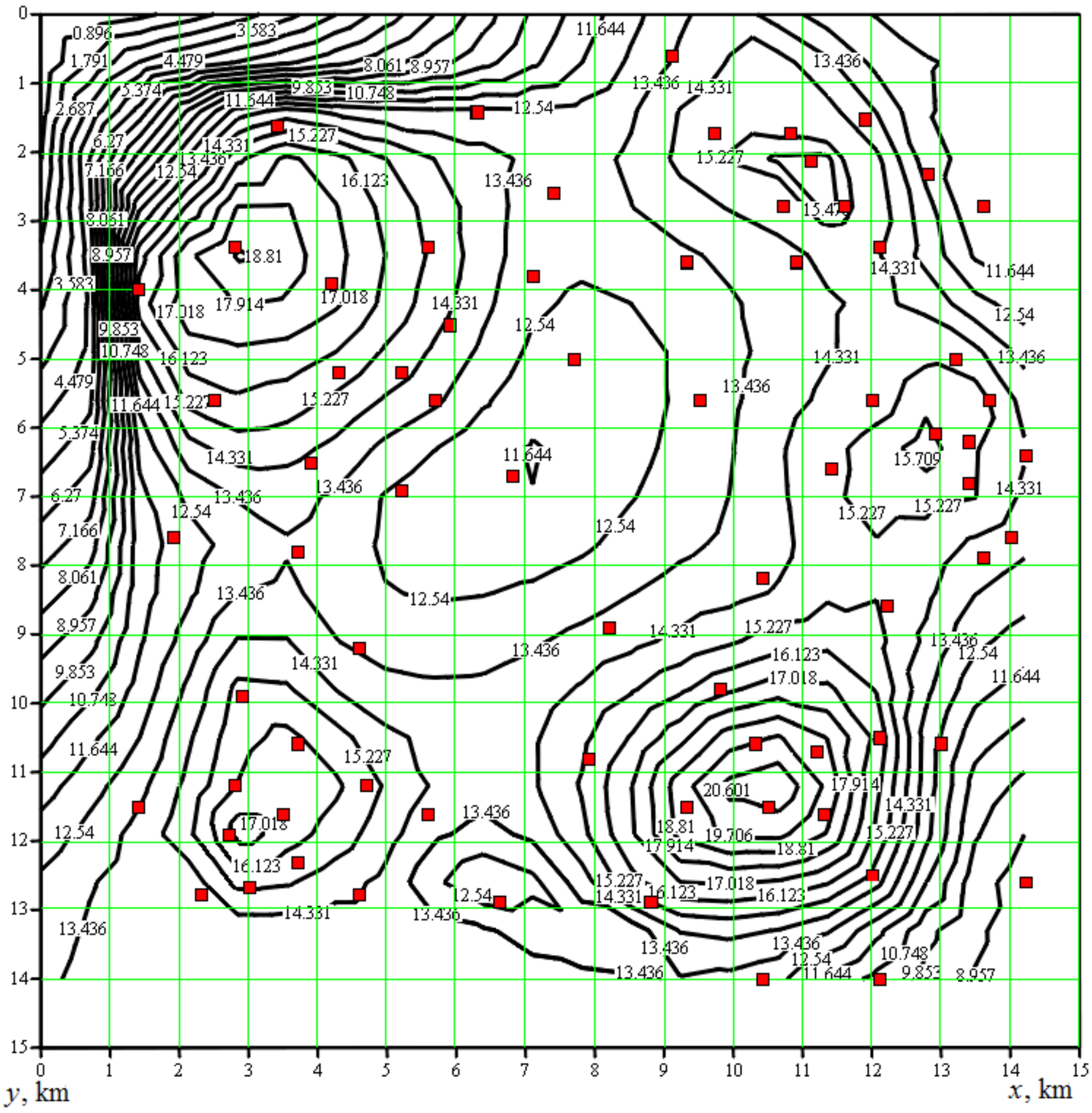

**Figure 10.** Isolines of intensity $\widetilde{V}_z(x, y, 0)$, mGal from several bodies constructed by $N = 73$ points of measurements (example 2). (Online version in colour.)

The regularization parameter is chosen to be $\alpha = 10^{-8}$, but its value influences weakly on the solution because the constraints on the solution play a role of regularization.

In table 3, the obtained solution is given.

**Table 3**. The obtained solution

| Body | $a$, km | ε | ρ, g/cm$^3$ | $x_0$, km | $y_0$, km | $z_0$, km | $\upsilon$, km$^3$ | $M$, bln t |
|---|---|---|---|---|---|---|---|---|
| 1 | 2.341 | 0.54 | 1.64 | 2.84 | 3.34 | 4.40 | 29.02 | 47.60 |
| 2 | 1.306 | 1.60 | 2.34 | 10.38 | 12.00 | 4.00 | 14.91 | 34.90 |
| 3 | 1.696 | 1.04 | 1.54 | 2.84 | 12.00 | 4.20 | 21.23 | 32.70 |
| 4 | 1.235 | 1.44 | 2.74 | 11.00 | 1.24 | 4.60 | 11.35 | 31.10 |
| 5 | 1.320 | 0.74 | 3.34 | 13.70 | 6.50 | 4.10 | 7.13 | 23.80 |

On the parameters of five bodies-spheroids adduced in table 3, the $z$-intensities $V_z(x, y, 0)$, $x \in [0,15]$ km, $y \in [0,15]$ km are calculated and the isolines of $V_z$ are constructed (figure 11).

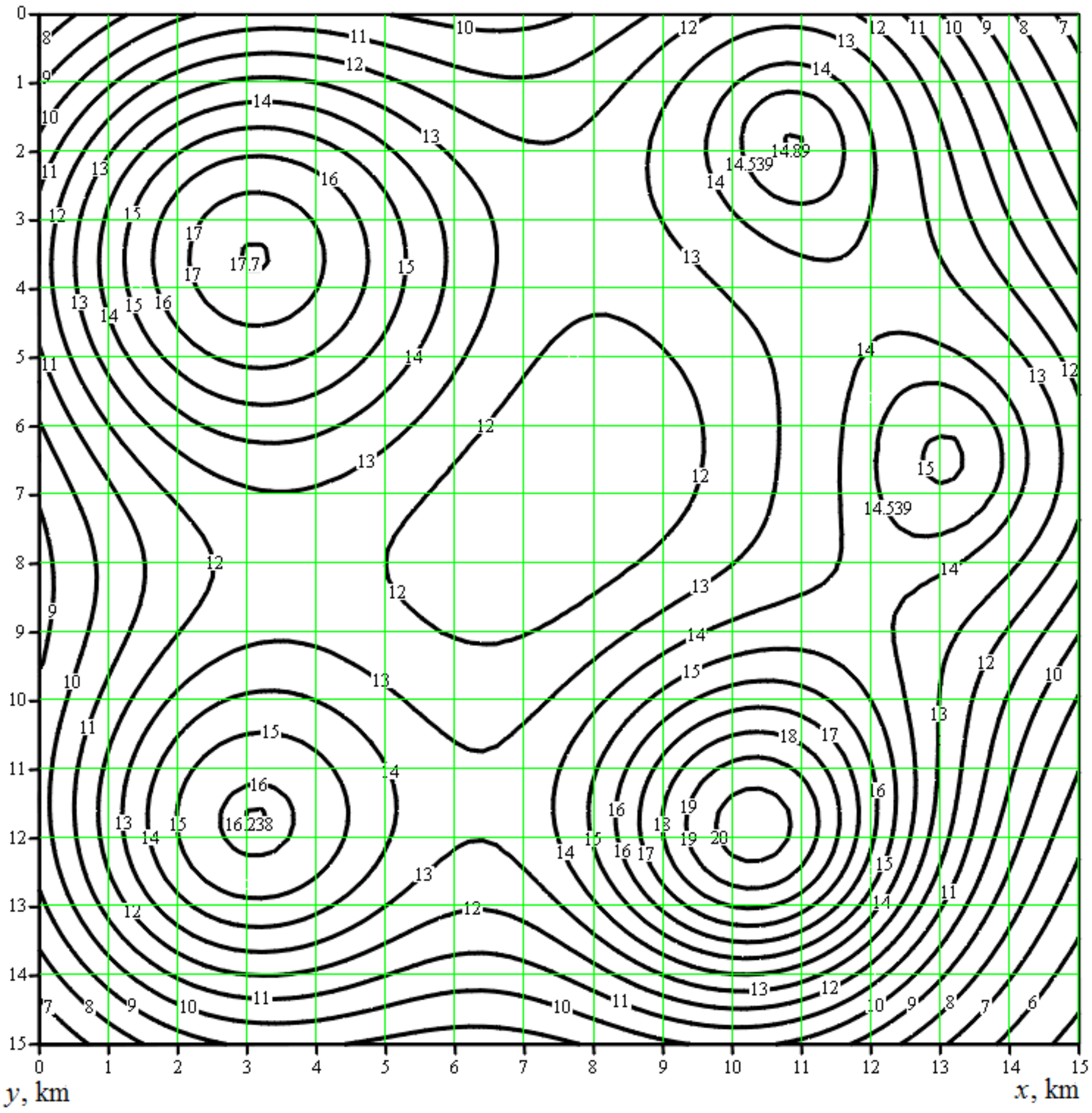

**Figure 11.** Isolines of intensity $V_z(x, y, 0)$, mGal from five spheroids (example 2). (Online version in colour.)

The comparison of isolines measured $\widetilde{V}_z$ (figure 10) and calculated $V_z$ (figure 11) shows that the main goal of modeling is reached, namely, five bodies are selected; their surface coordinates $x_0, y_0$ are estimated; the values of $V_z$ in five calculated 'poles' (figure 11) correspond to the values of $\widetilde{V}_z$ in measured 'poles' (figure 10), etc. in spite of the fact that example 2 is enough complicated for processing.

In figure 12, the images of five calculated spheroids are given. They can be identified with a deposit consisting of several bodies – ore, oil, intrusion, quartzites, shales, granite, basalt, et al.

## 7. Conclusion

Results obtained in the present paper and in [6–8] allow to draw the following conclusions:

1. Constraints on the solution (parameters of the deposit model), e.g., in the form of inequalities (5.3) help to eliminate nonuniqueness of the inverse gravimetry problem and enhance its stability. Use of the Tikhonov regularization method also promotes this (see (5.1), (5.2)).

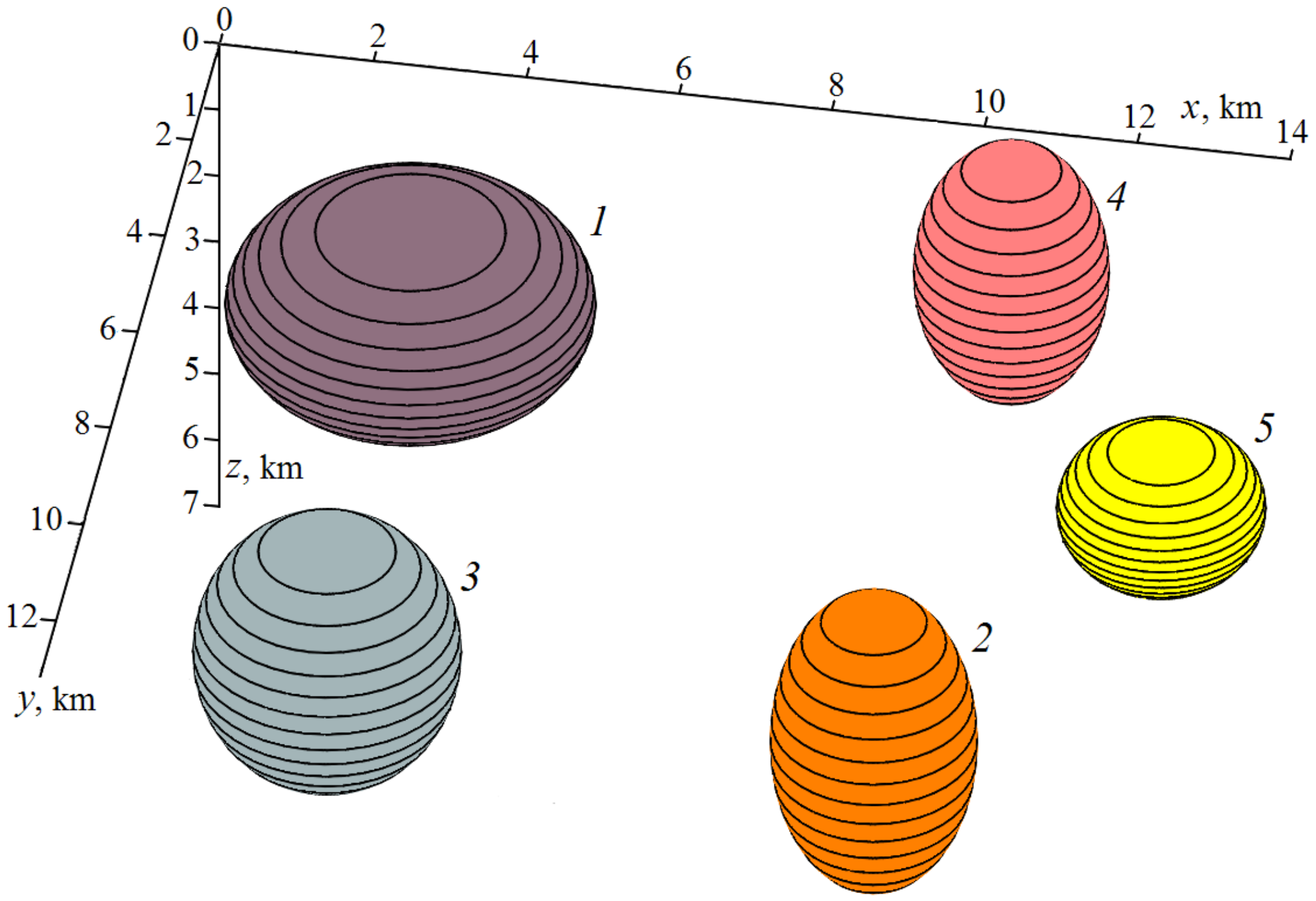

**Figure 12.** Five spheroids imitating the deposit in example 2. (Online version in colour.)

2. The solution of the inverse gravimetry problem may be stable (and unique) without regularization, but with use of rigid constraints on the solution.

3. The computerized realization of the problem of minimization of the functional $F_1$ or $F_2$ (see (5.1) or (5.2)) requires only one-dimensional arrays ($V_{z\,i}$, et al.) and not matrices (cf. [34]). The time for solving the problem of functional minimization by the coordinate descent method in the above-considered examples is of order 1–10 s at a CPU clock rate ~ 1 GHz. Thus, the proposed technique does not require large computer resources (memory and time).

4. General conclusion: the use of spheroids for deposit modeling, as well as the method for minimizing the smoothing functional with constraints on the sought parameters allows to solve the inverse gravimetry problem enough exactly, stably, uniquely and with small expenditures of computer memory and time even in the case of significant errors in the initial data.

This work was carried out with the support of the Russian Foundation for Basic Research, projects no. 09-08-00034 and 13-08-00442.

## References

1. Starostenko VI. 1978 *Stable numerical methods in gravimetry problems*. Kiev, Ukraine: Naukova Dumka.
2. Bulakh EG, Shinshin IV. 2000 Algorithmic and programmed solving the problem of constructing analytical model of gravitational field. *Geophys. J.* (Kiev) **22**, 107–114.
3. Pilkington M. 2014 Evaluating the utility of gravity gradient tensor components. *Geophysics* **79**, G1–G14. (doi:10.1190/geo2013-0130.1)
4. Salem A *et al.* 2014 Inversion of gravity data with isostatic constraints. *Geophysics* **79**, A45–A50. (doi:10.1190/geo2014-0220.1)
5. Golov IN, Sizikov VS. 2001 Inverse gravimetry problem as the intravising problem. *Sci-tech. Vestnik SPb IFMO* (*TU*). No. 3 (197), 171–175.
6. Golov IN, Sizikov VS. 2005 Modeling the inverse gravimetry problem with use of spheroids, nonlinear programming and regularization. *Geophys. J.* (Kiev) **27**, 454–462.

7. Golov IN, Sizikov VS. 2005 On correct solving the inverse gravimetry problem. *Rus. J. Geophysics* (St. Petersburg). No. 39–40, 84–91.
8. Golov IN., Sizikov VS. 2009 Modeling of deposits by spheroids // *Izv. Phys. Solid Earth* **45**, 258–271. (doi:10.1134/S1069351309030070)
9. Sizikov VS. 1967 Mass distribution in galaxies deduced from radial-velocity and photometry data. *Astrophysics* **3**, 124–128.
10. Vorontsov-Vel'yaminov BA. 1987 *Extragalactic astronomy*. Chur, Switzerland: Harwood Academic Publishers.
11. Bosold A *et al*. 2005 The structural geology of the High Central Zagros revisited (Iran). *Petroleum Geoscience* **11**, 225–238. (doi:10.1144/1354-079304-646)
12. Michel V, Wolf K. 2008 Numerical aspects of a spline-based multiresolution recovery of the harmonic mass density out of gravity functionals. *Geophys. J. International* **173**, 1–16. (doi:10.1111/j.1365-246X.2007.03700.x)
13. Salem A *et al*. 2013 Moho depth and sediment thickness estimation beneath the Red Sea derived from satellite and terrestrial gravity data. *Geophysics* **78**, G89–G101. (doi:10.1190/geo2012-0150.1)
14. Lavrent'ev MM, Romanov VG, Shishatskiĭ SP. 1986 *Ill-posed problems of mathematical physics and analysis*. Providence, USA: AMS.
15. Tikhonov AN, Arsenin VYa. 1977 *Solution of ill-posed problems*. New York, NY: Wiley.
16. Verlan' AF, Sizikov VS. 1986 *Integral equations: methods, algorithms, programs*. Kiev, Ukraine: Naukova Dumka.
17. Leonov AS. 2010 *Solving ill-posed inverse problems*. Moscow, Russia: Librokom.
18. Bakushinsky A, Goncharsky A. 1994 *Ill-posed problems: theory and applications*. Dordrecht, Germany: Kluwer.
19. Subbotin MF. 1949 *Course of celestial mechanics*, vol. 3. Leningrad-Moscow, Russia: GITTL.
20. Duboshin GN. 1975 *Celestial mechanics. Basic problems and methods*, 3rd edn. Moscow, Russia: Nauka.
21. Yun'kov AA, Afanas'ev NL. 1952 The direct and inverse problem of $\Delta g$ for a triaxial ellipsoid and an ellipsoid of revolution. *Izv. Dnepropetrovsk. Gorn. Inst.* **22**, 17–21 (Moscow-Khar'kov, SSSR: Ugletekhizdat).
22. Stepanova IE. 2001 A robust algorithm for reconstructing an ellipsoid. *Izv. Phys. Solid Earth* **37**, 955–960.
23. Krizskii VN, Gerasimov IA, Viktorov SV. 2002 Mathematical modeling inverse problems of potential geoelectric fields in axially symmetric piecewise homogeneous media. *Vestnik Zaporozh. Gosud. Univers. Fiz.-Mat. Nauki*. No. 1, 1–5.
24. *The physics of medical imaging*. 1988, in 2 vol. / ed. by S. Webb. Bristol, UK: IOP Publishing.
25. Sizikov VS. 2001 *Mathemathical methods for processing the results of measurements*. St. Petersburg, Russia: Politekhnika.
26. Marusina MYa, Kaznacheeva AO. 2007 Contemporary state and perspectives of development of tomography. *Sci-tech. Vestnik SPbSU IFMO*. No. 8 (42), 3–13.
27. Sizikov VS. 2011 *Inverse applied problems and MatLab*. St. Petersburg, Russia: Lan'.
28. MacLennan K, Karaoulis M, Revil A. 2014 Complex conductivity tomography using low-frequency crosswell electromagnetic data. *Geophysics* **79**, E23–E38. (doi:10.1190/geo 2012-0531.1)
29. Cho ZH, Jones JP, Singh M. 1993 *Foundations of medical imaging*. New York, USA: Wiley.
30. Sretensky LN. 1954 On the uniqueness of the determination of the shape of an attractive body from its outer potential. *Doklady Akademii Nauk SSSR* **99**, 21–22.

31. Bulakh EG, Markova MN. 2008 Inverse gravity problems for models composed of bodies of Sretenskii's class // *Izv. Phys. Solid Earth* **45**, 258–271. (doi:10.1134/S1069351308070033)
32. Bulakh EG, Rzhanitsyn VA, Markova MN. 1976 *Application of a minimization method for solving problems of structural geology by data of gravitational prospecting*. Kiev, Ukraine: Naukova Dumka.
33. Starostenko VI, Pashko VF, Zavorot'ko AN. 1992 Experience in solving strongly unstable inverse linear gravity problem. *Izv. Phys. Solid Earth*. No. 8, 24–44.
34. Strakhov VN. 2004 New paradigm in the theory of linear ill-posed problems adequate to requirements of geophysical practice. I–V. *Geophys. J.* (Kiev) **26**. No. 1–4.
35. Bertete-Aguirre H, Cherkaev E, Oristaglio M. 2002 Non-smooth gravity problem with total variation penalization functional. *Geophys. J. International* **149**, 499–507.
36. Zhuravlev IA. 1998 On solving the inverse gravimetry problem in the class of distributed density. *Geophys. J.* (Kiev) **20**, 86–94.
37. Himmelblau DM. 1972 *Applied nonlinear programming*. New York, USA: McGraw-Hill.
38. Vasiliev FP. 1981 *Methods for solving extreme problems*. Moscow, Russia: Nauka.
39. Mottl J, Mottlova L. 1984 The simultaneous solution of the inverse problem of gravimetry and magnetics by means of non-linear programming. *Geophis. J. RAS* **76**, 563–579.
40. D'yakonov VP. 1984 *Handbook on algorithms and programs in Basic language for PC*. Moscow, Russia: Nauka.